\documentclass[12pt]{article}
\usepackage{amsthm}
\usepackage{amsmath}
\usepackage{amssymb}
\usepackage{isomath}
\usepackage{graphicx}
\usepackage{float}
\usepackage[margin=1in]{geometry}
\usepackage{biblatex}
\title{The One-Seventh Area Triangle in the Complex Plane: Proving the Proportionality Using Features of the Complex Number System}
\author{Mathew Miltonhardy}
\date{September 2025}

\newcommand{\vertex}[1]{\boldsymbol{#1}}
\newcommand{\abs}[1]{\left\lvert #1 \right\rvert}
\let\oldsqrt\sqrt
\def\sqrt{\mathpalette\DHLhksqrt}
\def\DHLhksqrt#1#2{%
\setbox0=\hbox{$#1\oldsqrt{#2\,}$}\dimen0=\ht0
\advance\dimen0-0.2\ht0
\setbox2=\hbox{\vrule height\ht0 depth -\dimen0}%
{\box0\lower0.4pt\box2}}

\begin{document}

\maketitle

\section{Introduction}

This paper aims to provide an easy to follow proof that the area of the inner triangle formed by three cevians one-third along the edge of a given outer triangle is one-seventh the area of the outer triangle. This is a widely known problem, also called Feynman's triangle. There are many geometric proofs available online, but here we will use a purely algebraic approach. We will be using a triangle expressly in the complex plane, and use features of the complex plane to determine the vertices of the inner triangle formed. 

\begin{figure}[htb]
\makebox[\textwidth][c]{\includegraphics[width=\textwidth]{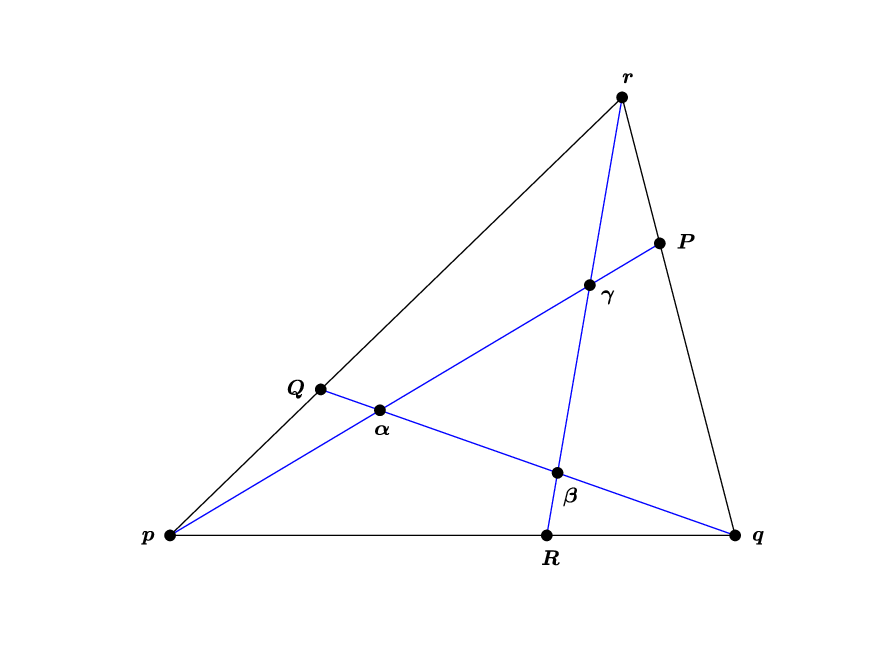}}
\caption{Triangle formed by vertices in the complex plane where $\boldsymbol{p}$ is at the origin, $\boldsymbol{q}$ is purely real, and $\boldsymbol{r}$ is any point in the complex plane. Cevians are drawn from each vertex to the point one-third along the opposite side in a clockwise direction.}
\label{fig:triangle-diced}
\end{figure}

To simplify this problem, we will be observing a triangle where one vertex at is the origin, one is purely real, and the other is anywhere in the complex plane, as seen in Figure \ref{fig:triangle-diced}. Once we have proven the proportionality true for this specific triangle, we will know it holds true for all triangles in the complex plane. 

This is because the area of a triangle is invariant under translations and rotations in the plane \cite{yr8Textbook}. Intuitively, if we were to take a piece of paper on a table, then slide it across a few centimetres and spin it a few degrees, the area of that paper would not change at all. The same can be done with a triangle in the complex plane. By translating and rotating our simple triangle, the area will not change.

To complete this proof, we must first find the area of the outer triangle, then the points in the complex plane of the vertices forming the inner triangle, then the area of the inner triangle. Hopefully after those steps, we will find that the area of the inner triangle is indeed one-seventh the area of the outer triangle.

\section{Area of the Outer Triangle}

\begin{figure}[htb]
\makebox[\textwidth][c]{\includegraphics[width=\textwidth]{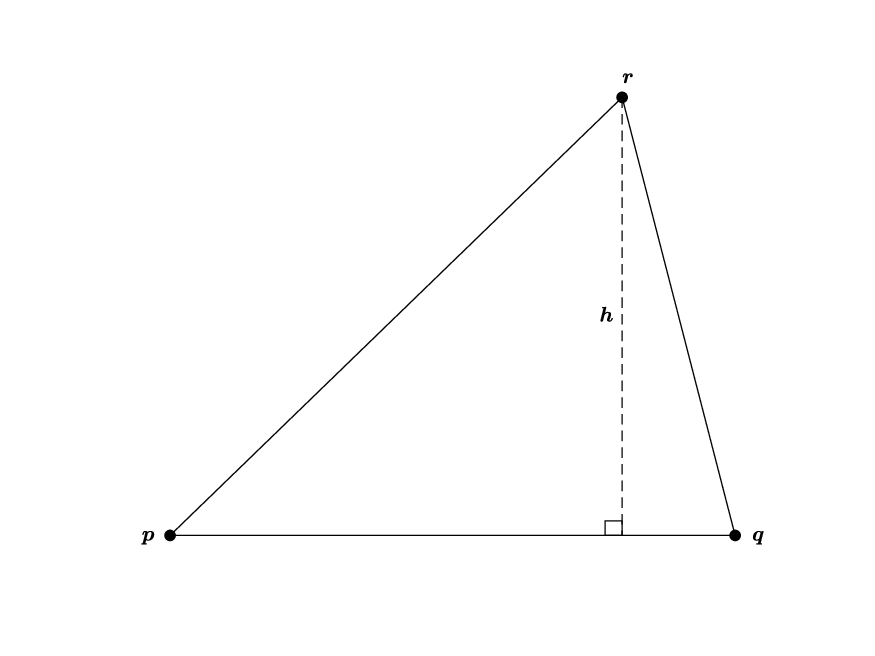}}
\caption{Triangle formed by vertices in the complex plane where $\vertex{p}$ is at the origin, $\vertex{q}$ is purely real, and $\vertex{r}$ is any point in the complex plane. The height of the triangle is labelled $\vertex{h}$.}
\label{fig:triangle-whole}
\end{figure}

Finding the area of the whole simple triangle is relatively straight forward. Since the edge of the triangle connecting vertices $\vertex{p}$ and $\vertex{q}$ is horizontal, we can find the height of the triangle as just the imaginary component of the vertex $\vertex{r}$. We can defined this height as $\vertex{h}$ shown in Figure \ref{fig:triangle-whole}. We will define the vertices as $\vertex{p} = 0 + 0i$, $\vertex{q} = b + 0i$, and $\vertex{r} = x + iy$.

Then we can calculate the area of the triangle as follows:

\begin{align*}
    A_{pqr} & = \frac{1}{2} \abs{\vertex{q} - \vertex{p}}  h \\
    & = \frac{1}{2} by
\end{align*}

While this is useful, it would be more useful to determine the area of any triangle formed in the complex plane. That is, a triangle where it's three non-collinear points lie in the complex plane.

\begin{figure}[htb]
\makebox[\textwidth][c]{\includegraphics[width=\textwidth]{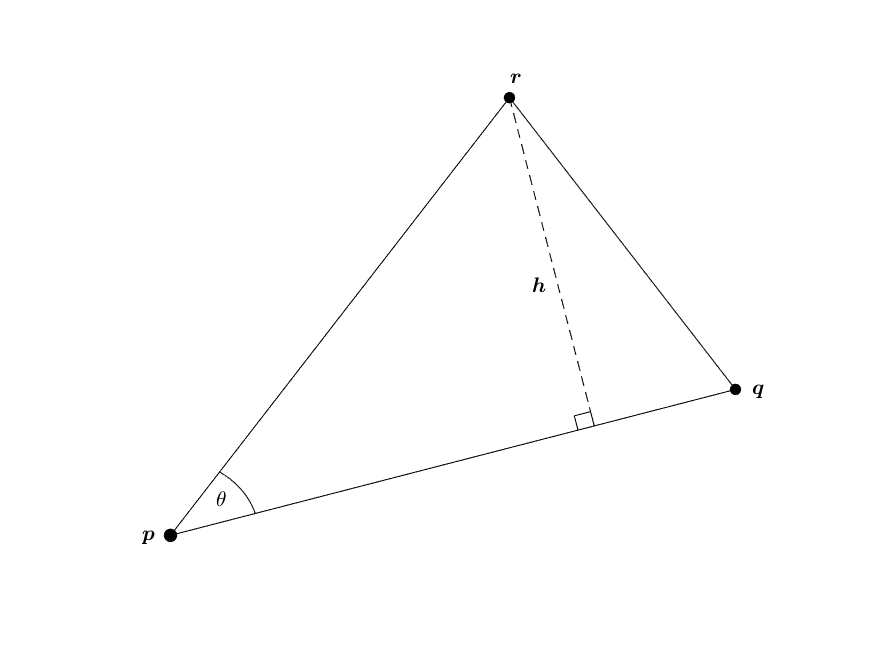}}
\caption{Triangle formed by vertices in the complex plane where $\vertex{p}$, $\vertex{q}$, and $\vertex{r}$ are any points in the complex plane. The height of the triangle is labelled $\vertex{h}$, and the angle $\angle \vertex{p}\vertex{q}\vertex{r}$ is labelled $\theta$.}
\label{fig:triangle-wholeRotated}
\end{figure}

To do this, we form a triangle with vertices $\vertex{p} = c + id$, $\vertex{q} = a + ib$, and $\vertex{r} = x + iy$, label the height of the triangle as $\vertex{h}$, and label the angle $\angle \vertex{p}\vertex{q}\vertex{r}$ as $\theta$, as shown in Figure \ref{fig:triangle-wholeRotated}.

To determine the value of $\vertex{h}$, we need to find the value of $\theta$. To do this, we notice that this angle will be the difference in argument of the vectors $\left(\vertex{r} - \vertex{p}\right)$ and $\left(\vertex{q} - \vertex{p}\right)$. Using this knowledge, we can find an equation for $\theta$ in terms of the values of our vertices such that,

\begin{align*}
    \theta &= \text{Arg}\left(\vertex{r} - \vertex{p}\right) - \text{Arg}\left(\vertex{q} - \vertex{p}\right) \\
    &= \text{Arg}\left(\frac{\vertex{r} - \vertex{p}}{\vertex{q} - \vertex{p}}\right) \\
    &= \text{Arg}\left(\frac{\left(x - c\right) + i \left(y - d\right)}{\left(a - c\right) + i \left(b-d\right)}\right) \\
    &= \text{Arg}\left(\frac{\left(x-c\right)\left(a-c\right) + \left(b-d\right)\left(y-d\right)}{\left(a-c\right)^2 + \left(b-d\right)^2} + i\frac{\left(a-c\right)\left(y-d\right) - \left(b-d\right)\left(x-c)\right)}{\left(a-c\right)^2 + \left(b-d\right)^2}\right) \\
    &= \arctan \left(\frac{\left(a-c\right)\left(y-d\right) - \left(b-d\right)\left(x-c\right)}{\left(x-c\right)\left(a-c\right) + \left(b-d\right)\left(y-d\right)}\right).
\end{align*}

Now that we know the value of $\theta$, we can use this to determine the value of $\vertex{h}$. As $\vertex{h}$ is perpendicular to the vector $\vertex{q} - \vertex{p}$, and thus forms a right-angled triangle, we can solve for $\vertex{h}$ as follows.

\begin{align*}
    \vertex{h} &= \abs{\vertex{r} - \vertex{p}}\sin{\theta} \\
    &= \abs{\vertex{r} - \vertex{p}} \sin \left(\arctan \left(\frac{\left(a-c\right)\left(y-d\right) - \left(b-d\right)\left(x-c\right)}{\left(x-c\right)\left(a-c\right) + \left(b-d\right)\left(y-d\right)}\right)\right) \\
    &= \abs{\vertex{r} - \vertex{p}} \cdot \frac{\left(a-c\right)\left(y-d\right) - \left(b-d\right)\left(x-c\right)}{\sqrt{\left[\left(a-c\right)\left(y-d\right) - \left(b-d\right)\left(x-c\right)\right]^2 + \left[\left(x-c\right)\left(a-c\right) + \left(b-d\right)\left(y-d\right)\right]^2 }} \\ 
    &= \abs{\vertex{r} - \vertex{p}} \cdot \frac{\left(a-c\right)\left(y-d\right) - \left(b-d\right)\left(x-c\right)}{\sqrt{\left(a-c\right)^2\left(y-d\right)^2 + \left(b-d\right)^2\left(x-c\right)^2 + \left(x-c\right)^2\left(a-c\right)^2 + \left(b-d\right)^2\left(y-d\right)^2}} \\
    &= \sqrt{\left(x-c\right)^2 + \left(y-d\right)^2} \cdot \frac{\left(a-c\right)\left(y-d\right) - \left(b-d\right)\left(x-c\right)}{\sqrt{\left(a-c\right)^2 + \left(b-d\right)^2} \cdot \sqrt{\left(x-c\right)^2 + \left(y-d\right)^2}} \\
    &= \frac{\left(a-c\right)\left(y-d\right) - \left(b-d\right)\left(x-c\right)}{\sqrt{\left(a-c\right)^2 + \left(b-d\right)^2}}.
\end{align*}

Now, a commonly known formula for the area of a triangle is one-half of the base multiplied by the height, and we now know the height, and know that the base of the triangle is the magnitude of the vector $\vertex{q} - \vertex{p}$, that is, $\abs{\vertex{q} - \vertex{p}}$. This gives us the formula for area as

\begin{align*}
    A_{pqr} &= \frac{1}{2} \abs{\vertex{q} - \vertex{p}} \vertex{h} \\
    &= \frac{1}{2} \sqrt{\left(a-c\right)^2 + \left(b-d\right)^2} \cdot \frac{\left(a-c\right)\left(y-d\right) - \left(b-d\right)\left(x-c\right)}{\sqrt{\left(a-c\right)^2 + \left(b-d\right)^2}} \\
    &= \frac{1}{2}\left(\left(a-c\right)\left(y-d\right) - \left(b-d\right)\left(x-c\right)\right) \\
    &= \frac{1}{2}\left(ay - ad - cy + cd - bx + bc + dx - cd\right) \\
    &= \frac{1}{2}\left(c\left(b-y\right) - d\left(a-x\right) + \left(ay-bx\right)\right) \\
    &= \frac{1}{2} \begin{vmatrix}
        c & d & 1 \\
        a & b & 1 \\
        x & y & 1
    \end{vmatrix}.
\end{align*}

We can see that the equation for area is the determinant of a matrix formed by the real and imaginary parts of each of the vertices of the triangle. This is known as the area of a triangle in determinant form and is generalised as below, where $x_i$ and $y_i$ represent the real and imaginary parts, respectively, of the vertices of a triangle \cite{determinantFormula}.

\begin{equation}
    A = \frac{1}{2} \begin{vmatrix}
        x_1 & y_1 & 1 \\
        x_2 & y_2 & 1 \\
        x_3 & y_3 & 1
    \end{vmatrix}.
    \label{eqn:areaDetForm}
\end{equation}

\section{Determining Coordinates}

To find the area of the inner triangle, we need to find the points in the complex plane of the vertices $\vertex{\alpha}$, $\vertex{\beta}$, and $\vertex{\gamma}$. But to do this, we first need the find the points along the edges of the triangle, $\vertex{P}$, $\vertex{Q}$, and $\vertex{R}$. To do this, we will revert to using our simple triangle with vertices $\vertex{p} = 0+0i$, $\vertex{q} = b + 0i$, and $\vertex{r} = x+iy$.

\subsection{The Edges}

To find the point $\vertex{P}$, we start at the vertex $\vertex{r}$, then move one third along the edge towards $\vertex{q}$. That is, we move one-third of the vector $\vertex{q} - \vertex{r}$. This shows us that

\begin{align*}
    \vertex{P} & = \vertex{r} + \frac{1}{3} \left(\vertex{q} - \vertex{r}\right) \\
    & = x + iy + \frac{1}{3} \left(b - x - iy\right) \\
    & = \frac{2}{3}x + \frac{1}{3}b + \frac{2}{3}yi.
\end{align*}

Using the same reasoning, we can see that

\begin{align*}
    \vertex{Q} & = \vertex{p} + \frac{1}{3} \left(\vertex{r} - \vertex{p}\right) \\
    & = \frac{1}{3} \left(x + iy\right) \\
    & = \frac{1}{3}x + \frac{1}{3}yi.
\end{align*}

To find the point $\vertex{R}$, we will instead start at vertex $\vertex{p}$ and moved two-thirds along the edge towards $\vertex{q}$, that is, two-thirds along the vector $\vertex{q} - \vertex{p}$. This shows us that

\begin{align*}
    \vertex{R} & = \vertex{p} + \frac{2}{3} \left(\vertex{q} - \vertex{p}\right) \\
    & = \frac{2}{3} \left(b + 0i\right) \\
    & = \frac{2}{3}b.
\end{align*}

\subsection{The Inner Vertices}

To determine the inner vertices, we must solve a set of two simultaneous equations for each vertex. We will use a trick to make this process easier for us and eliminate a portion of the calculations.

\subsubsection{Determining $\alpha$}

We can see that the vertex $\vertex{\alpha}$, can be written as a linear combination of the points and vectors we have found so far, such that:

\begin{align*}
    \vertex{\alpha} & = \vertex{p} + u\left(\vertex{P} - \vertex{p}\right), \text{ and } \\
    \vertex{\alpha} & = \vertex{q} + v\left(\vertex{Q} - \vertex{q}\right),
\end{align*}

where $u,v \in \mathbb{R}$.

We will solve this equation for $u$. To do this, we will equate both sides, and then multiply through by the complex conjugate of $\left(\vertex{Q} - \vertex{q}\right)$, that is, $\overline{\left(\vertex{Q} - \vertex{q}\right)}$. This yields

\begin{align}
    u\vertex{P}\overline{\left(\vertex{Q} - \vertex{q}\right)} = \vertex{q} \overline{\left(\vertex{Q} - \vertex{q}\right)} + v\left(\vertex{Q} - \vertex{q}\right)\overline{\left(\vertex{Q} - \vertex{q}\right)}.
    \label{eqn:findingAlpha}
\end{align}

Evaluating the left hand side of (\ref{eqn:findingAlpha}), we get

\begin{align*}
    LHS & = u\left(\frac{2}{3}x + \frac{1}{3}b + \frac{2}{3}yi\right)\left(\frac{1}{3}x - b - \frac{1}{3}yi\right) \\ 
    & = u\left(\left(\frac{2}{3}x + \frac{1}{3}b\right)\left(\frac{1}{3}x - b\right) + i\left(\frac{2}{9}xy - \frac{2}{3}by\right) - i\left(\frac{2}{9}xy + \frac{1}{9}by\right) + \frac{2}{9}y^2\right) \\
    & = u\left(\left(\frac{2}{3}x + \frac{1}{3}b\right)\left(\frac{1}{3}x - b\right) + \frac{2}{9}y^2 \right) + i\left(-\frac{7}{9}uby\right).
\end{align*}

We will leave that as it is for now, and will begin evaluating the right hand side, showing

\begin{align*}
    RHS & = b\left(\frac{1}{3}x - b - \frac{1}{3}yi\right) + v\left(\vertex{Q} - \vertex{q}\right)\overline{\left(\vertex{Q} - \vertex{q}\right)} \\
    & = b\left(\frac{1}{3}x - b\right) + v\abs{\vertex{Q} - \vertex{q}}^2 - \frac{1}{3}byi.
\end{align*}

Equating the imaginary parts of the left and right hand sides allows us to remove the unknown variable $v$, as $\abs{\vertex{Q} - \vertex{q}}^2$ is purely real. This gives us 

\begin{align*}
    \mathfrak{Im}(LHS) & = \mathfrak{Im}(RHS) \\
    -\frac{7}{9}uby & = -\frac{1}{3}by \\
    \frac{7}{3}u & = 1 \\
    u = \frac{3}{7}.
\end{align*}

Therefore, we now know that $\vertex{\alpha} = \frac{3}{7} \vertex{P}$.

\subsubsection{Determining $\beta$}

Using the same logic as above, we can see that $\vertex{\beta}$ can be written as a linear combination of the points and vectors found so far, such that:

\begin{align*}
    \vertex{\beta} & = \vertex{q} + u\left(\vertex{Q} - \vertex{q}\right), \text{ and } \\
    \vertex{\beta} & = \vertex{r} + v\left(\vertex{R} - \vertex{r}\right),
\end{align*}

where $u,v \in \mathbb{R}$.

We will again solve this equation for $u$ by equating both sides and multiplying through by a complex conjugate, this time using $\overline{\left(\vertex{R} - \vertex{r}\right)}$. This yields

\begin{align}
    \vertex{q} \overline{\left(\vertex{R} - \vertex{r}\right)} + u\left(\vertex{Q} - \vertex{q}\right)\overline{\left(\vertex{R} - \vertex{r}\right)} = \vertex{r}\overline{\left(\vertex{R} - \vertex{r}\right)} + v\left(\vertex{R} - \vertex{r}\right)\overline{\left(\vertex{R} - \vertex{r}\right)}.
    \label{eqn:findingBeta}
\end{align}

Evaluating the left have side of \ref{eqn:findingBeta} we get 

\begin{align*}
    LHS & = b \left(\frac{2}{3}b - x + iy\right) + u\left(\frac{1}{3}x - b + \frac{1}{3}yi\right)\left(\frac{2}{3}b - x + iy\right) \\
    & = b\left(\frac{2}{3}b - x\right) + iby + u\left(\left(\frac{1}{3}x - b\right)\left(\frac{2}{3}b - x\right) + i\left(\frac{2}{9}by - \frac{1}{3}xy\right) + i\left(\frac{1}{3}xy - by\right) - \frac{1}{3}y^2\right) \\
    & = b\left(\frac{2}{3}b - x\right) + u\left(\left(\frac{1}{3}x - b\right)\left(\frac{2}{3}b - x\right) - \frac{1}{3}y^2 \right) + i\left(by + u\left(- \frac{7}{9}by\right)\right).
\end{align*}

We will leave that as it is, since we have the imaginary part of the equation in a simple form. Now evaluating the right hand side we get

\begin{align*}
    RHS & = \left(x + iy\right)\left(\frac{2}{3}b - x + iy\right) + v\left(\vertex{R} - \vertex{r}\right)\overline{\left(\vertex{R} - \vertex{r}\right)} \\
    & = x\left(\frac{2}{3}b - x\right) +i\left(\frac{2}{3}by - xy\right) + ixy - y^2 + v\abs{\vertex{R} - \vertex{r}}^2 \\
    & = x\left(\frac{2}{3}b - x\right) - y^2 + v\abs{\vertex{R} - \vertex{r}}^2 + i\left(\frac{2}{3}by\right).
\end{align*}

Equating the imaginary parts of the left and right hand sides allows us to remove one of the unknown variables, $v$, as $\abs{\vertex{R} - \vertex{r}}^2$ is purely real. This gives us

\begin{align*}
    \mathfrak{Im}(LHS) & = \mathfrak{Im}(RHS) \\
    by + u\left(-\frac{7}{9}by\right) & = \frac{2}{3}by \\
    \frac{7}{3}u\left(-\frac{1}{3}by\right) & = -\frac{1}{3}by \\
    \frac{7}{3}u & = 1 \\
    u & = \frac{3}{7}.
\end{align*}

Therefore, we now know that $\vertex{\beta} = \vertex{q} + \frac{3}{7}\left(\vertex{Q} - \vertex{q}\right)$.

\subsubsection{Determining $\gamma$}

Using the same logic yet again, we can write $\vertex{\gamma}$ as a linear combination of the points and vectors found so far such that:

\begin{align*}
    \vertex{\gamma} & = \vertex{p} + u\left(\vertex{P} - \vertex{p}\right), \text{ and } \\
    \vertex{\gamma} & = \vertex{r} + v\left(\vertex{R} - \vertex{r}\right),
\end{align*}

where $u,v \in \mathbb{R}$.

By equating both sides of the equation and multiplying through by the complex conjugate $\overline{\left(\vertex{R} - \vertex{r}\right)}$, we can solve for $u$. This yields

\begin{align}
    u\vertex{P}\overline{\left(\vertex{R} - \vertex{r}\right)} = \vertex{r}\overline{\left(\vertex{R} - \vertex{r}\right)} + v\left(\vertex{R} - \vertex{r}\right)\overline{\left(\vertex{R} - \vertex{r}\right)}.
    \label{eqn:findingGamma}
\end{align}

Evaluating the left have side of \ref{eqn:findingGamma} we get

\begin{align*}
    LHS & = u\left(\frac{2}{3}x + \frac{1}{3}b + \frac{2}{3}yi\right)\left(\frac{2}{3}b - x + iy\right) \\ 
    & = u\left(\left(\frac{2}{3}x + \frac{1}{3}b\right)\left(\frac{2}{3}b - x\right) + i\left(\frac{4}{9}by - \frac{2}{3}xy\right) + i\left(\frac{2}{3}xy + \frac{1}{3}by\right) - \frac{2}{3}y^2\right) \\ 
    & = u\left(\left(\frac{2}{3}x + \frac{1}{3}b\right)\left(\frac{2}{3}b - x\right) - \frac{2}{3}y^2 \right) + i\left(\frac{7}{9}uby\right).
\end{align*}

We can leave that as it is because we have the imaginary part of the equation in a simple form. Now we will evaluate the right hand side to show

\begin{align*}
    RHS & = \left(x + iy\right)\left(\frac{2}{3}b - x + iy\right) + v\left(\vertex{R} - \vertex{r}\right)\overline{\left(\vertex{R} - \vertex{r}\right)} \\
    & = x\left(\frac{2}{3}b - x\right) + i\left(\frac{2}{3}by - xy\right) + ixy - y^2 + v\abs{\vertex{R} - \vertex{r}}^2 \\
    & = x\left(\frac{2}{3}b - x\right) - y^2 + v\abs{\vertex{R} - \vertex{r}}^2 + i\left(\frac{2}{3}by\right).
\end{align*}

As done previously, we equate the imaginary parts of both the left and right hand sides which allows us to remove the unknown variable $v$, as $\abs{\vertex{R} - \vertex{r}}^2$ is purely real. This gives us

\begin{align*}
    \mathfrak{Im}(LHS) & = \mathfrak{Im}(RHS) \\
    \frac{7}{9}uby & = \frac{2}{3}by \\
    \frac{7}{3}u & = 2 \\
    u & = \frac{6}{7}.
\end{align*}

Therefore, we now know that $\vertex{\gamma} = \frac{6}{7}\vertex{P}$. This means we now know all of the inner vertices and can use these to determine the area of the inner triangle.

\section{Area of the Inner Triangle}

To make the following calculations easier, we will rewrite our vertices and points in the form $z = \left(\mathfrak{Re}\left(z\right), \mathfrak{Im}\left(z\right)\right)$. This gives all of our vertices and points as below.

\begin{alignat*}{3}
\vertex{p} &= \left(0, \;0\right)  &\qquad\;\;  \vertex{P} &= \left(\frac{2}{3}x + \frac{1}{3}b, \; \frac{2}{3}y\right)  &\qquad\;\;  \vertex{\alpha} &= \frac{3}{7}\vertex{P} = \left(\frac{2}{7}x + \frac{1}{7}b, \; \frac{2}{7}y\right) \\
\vertex{q} &= \left(b, \; 0 \right)  &  \vertex{Q} &= \left(\frac{1}{3}x, \; \frac{1}{3}y \right) &  \vertex{\beta} &= \vertex{q} + \frac{3}{7}\left(\vertex{Q} - \vertex{q}\right) = \left(\frac{1}{7}x + \frac{4}{7}b, \; \frac{1}{7}y \right)\\
\vertex{r} &= \left(x, \; y\right)  &  \vertex{R} &= \left(\frac{2}{3}b, \; 0\right)  &  \vertex{\gamma} &= \frac{6}{7}\vertex{P} = \left(\frac{4}{7}x + \frac{2}{7}b, \; \frac{4}{7}y\right)
\end{alignat*}

Since we have the vertices of the the inner triangle, we can use (\ref{eqn:areaDetForm}) derived earlier to find the area of our triangle.

First, to verify our previous derivation, we will calculate the area of the whole triangle using (\ref{eqn:areaDetForm}). This gives us

\begin{align*}
    A_{pqr} &= \frac{1}{2} 
    \begin{vmatrix}
        0 & 0 & 1 \\
        b & 0 & 1 \\
        x & y & 1
    \end{vmatrix} \\
    &= \frac{1}{2}\left(by - 0\right) \\
    &= \frac{1}{2}by,
\end{align*}

which matches what we previously obtained.

Now, we will use (\ref{eqn:areaDetForm}) to find the area of our inner triangle. This gives us

\begin{align*}
    A_{\alpha \beta \gamma} &= \frac{1}{2} 
    \begin{vmatrix}
        \frac{2}{7}x + \frac{1}{7}b & \frac{2}{7}y & 1 \\
        \frac{1}{7}x + \frac{4}{7}b & \frac{1}{7}y & 1 \\
        \frac{4}{7}x + \frac{2}{7}b & \frac{4}{7}y & 1
    \end{vmatrix} \\
    &= \frac{1}{2}\left(\left(\frac{2}{7}x + \frac{1}{7}b\right)\left(\frac{1}{7}y - \frac{4}{7}y\right) - \frac{2}{7}y\left(\frac{1}{7}x + \frac{4}{7}b - \frac{4}{7}x - \frac{2}{7}b\right) + \frac{4}{7}y\left(\frac{1}{7}x + \frac{4}{7}b\right) - \frac{1}{7}y\left(\frac{4}{7}x + \frac{2}{7}b\right)\right) \\
    &= \frac{1}{2}\left(-\frac{3}{7}y\left(\frac{2}{7}x + \frac{1}{7}b\right) - \frac{2}{7}y\left(-\frac{3}{7}x + \frac{2}{7}b\right) + \frac{4}{7}y\left(\frac{1}{7}x + \frac{4}{7}b\right) - \frac{1}{7}y\left(\frac{4}{7}x + \frac{2}{7}b\right)\right) \\
    &= \frac{1}{7} \cdot \frac{1}{2}y\left(-\frac{6}{7}x - \frac{3}{7}b + \frac{6}{7}x - \frac{4}{7}b + \frac{4}{7}x +\frac{16}{7}b - \frac{4}{7}x - \frac{2}{7}b\right) \\
    &= \frac{1}{7} \cdot \frac{1}{2}yb\left(\frac{16-3-4-2}{7}\right) \\
    &= \frac{1}{7} \cdot \frac{1}{2}by \\
    &= \frac{1}{7} A_{pqr}.
\end{align*}

\section{Conclusion}

We have proven that the area of the inner triangle formed by three cevians one-third along the edge of the specific triangle in the complex plane where one vertex is the origin, one is purely real, and one is any complex number, is one-seventh the area of the outer triangle. We know that the area of a triangle is invariant under translations and rotations in the plane \cite{yr8Textbook}.

Any triangle in the complex plane with three non-collinear vertices can be translated such that one of the vertices lies at the origin, and then rotated about the origin such that one edge of the triangle lies on the real-axis. These translations and rotations do not change the area of the triangle, meaning that the area of the inner triangle will be one-seventh of the area of the outer triangle for any triangle in the complex plane.

We have derived a formula for the area of a triangle with any three non-collinear vertices in the complex plane, and as such, can find the area of any inner triangle formed by three cevians one-third along the edge of a given outer triangle. We know that the area of the inner triangle is simply one-seventh of the area of the outer triangle.

\section{Acknowledgements}

I would like to thank my second year university mathematics teacher, Frank V, who inspired this paper by including the one-seventh triangle problem in an assignment. I must also extend my gratitude to my cat, Inky, for his unwavering support and comfort throughout writing this paper.

\section{Further Research}

An interesting topic of further study would be to determine whether a similar proportionality can be derived when creating the cevians moving a different fraction along the edge of the triangle. Investigating initially what the proportionality would be when moving one-fourth along the edge of a given triangle, and then one-fifth, then one-sixth, and so on, would allow us to determine whether the proportion of area is constant, and if not, how to determine it. Would it be possible to derive a general solution for the proportion of area based on the fraction of the edge covered by the cevians? 

\printbibliography

\end{document}